\documentstyle[12pt,bezier]{article} 
\newcommand{\qbezier}{\bezier{120}}
\def\DottedCircle{
\bezier{4}(0.966,-0.259)(1.04,0)(0.966,0.259)
\bezier{4}(0.966,0.259)(0.897,0.518)(0.707,0.707)
\bezier{4}(0.707,0.707)(0.518,0.897)(0.259,0.966)
\bezier{4}(0.259,0.966)(0,1.04)(-0.259,0.966)
\bezier{4}(-0.259,0.966)(-0.518,0.897)(-0.707,0.707)
\bezier{4}(-0.707,0.707)(-0.897,0.518)(-0.966,0.259)
\bezier{4}(-0.966,0.259)(-1.04,0)(-0.966,-0.259)
\bezier{4}(-0.966,-0.259)(-0.897,-0.518)(-0.707,-0.707)
\bezier{4}(-0.707,-0.707)(-0.518,-0.897)(-0.259,-0.966)
\bezier{4}(-0.259,-0.966)(0,-1.04)(0.259,-0.966)
\bezier{4}(0.259,-0.966)(0.518,-0.897)(0.707,-0.707)
\bezier{4}(0.707,-0.707)(0.897,-0.518)(0.966,-0.259)
}
\def\Endpoint[#1]{
\ifcase#1
\put(1,0){\circle*{0.15}}
\or\put(0.866,0.5){\circle*{0.15}}
\or\put(0.5,0.866){\circle*{0.15}}
\or\put(0,1){\circle*{0.15}}
\or\put(-0.5,0.866){\circle*{0.15}}
\or\put(-0.866,0.5){\circle*{0.15}}
\or\put(-1,0){\circle*{0.15}}
\or\put(-0.866,-0.5){\circle*{0.15}}
\or\put(-0.5,-0.866){\circle*{0.15}}
\or\put(0,-1){\circle*{0.15}}
\or\put(0.5,-0.866){\circle*{0.15}}
\or\put(0.866,-0.5){\circle*{0.15}}
\fi}
\def\Arc[#1]{
\thicklines			
\ifcase#1
\bezier{25}(0.966,-0.259)(1.04,0)(0.966,0.259)
\or
\bezier{25}(0.966,0.259)(0.897,0.518)(0.707,0.707)
\or
\bezier{25}(0.707,0.707)(0.518,0.897)(0.259,0.966)
\or
\bezier{25}(0.259,0.966)(0,1.04)(-0.259,0.966)
\or
\bezier{25}(-0.259,0.966)(-0.518,0.897)(-0.707,0.707)
\or
\bezier{25}(-0.707,0.707)(-0.897,0.518)(-0.966,0.259)
\or
\bezier{25}(-0.966,0.259)(-1.04,0)(-0.966,-0.259)
\or
\bezier{25}(-0.966,-0.259)(-0.897,-0.518)(-0.707,-0.707)
\or
\bezier{25}(-0.707,-0.707)(-0.518,-0.897)(-0.259,-0.966)
\or
\bezier{25}(-0.259,-0.966)(0,-1.04)(0.259,-0.966)
\or
\bezier{25}(0.259,-0.966)(0.518,-0.897)(0.707,-0.707)
\or
\bezier{25}(0.707,-0.707)(0.897,-0.518)(0.966,-0.259)
\fi}
\def\DottedArc[#1]{
\ifcase#1
\bezier{4}(0.966,-0.259)(1.04,0)(0.966,0.259)
\or
\bezier{4}(0.966,0.259)(0.897,0.518)(0.707,0.707)
\or
\bezier{4}(0.707,0.707)(0.518,0.897)(0.259,0.966)
\or
\bezier{4}(0.259,0.966)(0,1.04)(-0.259,0.966)
\or
\bezier{4}(-0.259,0.966)(-0.518,0.897)(-0.707,0.707)
\or
\bezier{4}(-0.707,0.707)(-0.897,0.518)(-0.966,0.259)
\or
\bezier{4}(-0.966,0.259)(-1.04,0)(-0.966,-0.259)
\or
\bezier{4}(-0.966,-0.259)(-0.897,-0.518)(-0.707,-0.707)
\or
\bezier{4}(-0.707,-0.707)(-0.518,-0.897)(-0.259,-0.966)
\or
\bezier{4}(-0.259,-0.966)(0,-1.04)(0.259,-0.966)
\or
\bezier{4}(0.259,-0.966)(0.518,-0.897)(0.707,-0.707)
\or
\bezier{4}(0.707,-0.707)(0.897,-0.518)(0.966,-0.259)
\fi}
\def\Chord[#1,#2]{
\thinlines
\ifnum#1>#2\Chord[#2,#1]
\else\ifnum#1<#2
\ifcase#1
\ifcase#2
\or\qbezier(1,0)(0.516,0.138)(0.866,0.5)
\or\qbezier(1,0)(0.45,0.26)(0.5,0.866)
\or\qbezier(1,0)(0.327,0.327)(0,1)
\or\qbezier(1,0)(0.179,0.311)(-0.5,0.866)
\or\qbezier(1,0)(0.0536,0.2)(-0.866,0.5)
\or\put(1, 0){\line(-2, 0){2}}
\or\qbezier(1,0)(0.0536,-0.2)(-0.866,-0.5)
\or\qbezier(1,0)(0.179,-0.311)(-0.5,-0.866)
\or\qbezier(1,0)(0.327,-0.327)(0,-1)
\or\qbezier(1,0)(0.45,-0.26)(0.5,-0.866)
\or\qbezier(1,0)(0.516,-0.138)(0.866,-0.5)
\fi
\or\ifcase#2\or
\or\qbezier(0.866,0.5)(0.378,0.378)(0.5,0.866)
\or\qbezier(0.866,0.5)(0.26,0.45)(0,1)
\or\qbezier(0.866,0.5)(0.12,0.446)(-0.5,0.866)
\or\qbezier(0.866,0.5)(0,0.359)(-0.866,0.5)
\or\qbezier(0.866,0.5)(-0.0536,0.2)(-1,0)
\or\put(0.866, 0.5){\line(-5, -3){1.73}}
\or\qbezier(0.866,0.5)(0.146,-0.146)(-0.5,-0.866)
\or\qbezier(0.866,0.5)(0.311,-0.179)(0,-1)
\or\qbezier(0.866,0.5)(0.446,-0.12)(0.5,-0.866)
\or\qbezier(0.866,0.5)(0.52,0)(0.866,-0.5)
\fi
\or\ifcase#2\or\or
\or\qbezier(0.5,0.866)(0.138,0.516)(0,1)
\or\qbezier(0.5,0.866)(0,0.52)(-0.5,0.866)
\or\qbezier(0.5,0.866)(-0.12,0.446)(-0.866,0.5)
\or\qbezier(0.5,0.866)(-0.179,0.311)(-1,0)
\or\qbezier(0.5,0.866)(-0.146,0.146)(-0.866,-0.5)
\or\put(0.5, 0.866){\line(-3, -5){1}}
\or\qbezier(0.5,0.866)(0.2,-0.0536)(0,-1)
\or\qbezier(0.5,0.866)(0.359,0)(0.5,-0.866)
\or\qbezier(0.5,0.866)(0.446,0.12)(0.866,-0.5)
\fi
\or\ifcase#2\or\or\or
\or\qbezier(0,1.)(-0.138,0.516)(-0.5,0.866)
\or\qbezier(0,1.)(-0.26,0.45)(-0.866,0.5)
\or\qbezier(0,1.)(-0.327,0.327)(-1,0)
\or\qbezier(0,1.)(-0.311,0.179)(-0.866,-0.5)
\or\qbezier(0,1.)(-0.2,0.0536)(-0.5,-0.866)
\or\put(0, 1){\line(0, -2){2}}
\or\qbezier(0,1.)(0.2,0.0536)(0.5,-0.866)
\or\qbezier(0,1.)(0.311,0.179)(0.866,-0.5)
\fi
\or\ifcase#2\or\or\or\or
\or\qbezier(-0.5,0.866)(-0.378,0.378)(-0.866,0.5)
\or\qbezier(-0.5,0.866)(-0.45,0.26)(-1,0)
\or\qbezier(-0.5,0.866)(-0.446,0.12)(-0.866,-0.5)
\or\qbezier(-0.5,0.866)(-0.359,0)(-0.5,-0.866)
\or\qbezier(-0.5,0.866)(-0.2,-0.0536)(0,-1)
\or\put(-0.5, 0.866){\line(3, -5){1}}
\or\qbezier(-0.5,0.866)(0.146,0.146)(0.866,-0.5)
\fi
\or\ifcase#2\or\or\or\or\or
\or\qbezier(-0.866,0.5)(-0.516,0.138)(-1,0)
\or\qbezier(-0.866,0.5)(-0.52,0)(-0.866,-0.5)
\or\qbezier(-0.866,0.5)(-0.446,-0.12)(-0.5,-0.866)
\or\qbezier(-0.866,0.5)(-0.311,-0.179)(0,-1)
\or\qbezier(-0.866,0.5)(-0.146,-0.146)(0.5,-0.866)
\or\put(-0.866, 0.5){\line(5, -3){1.73}}
\fi
\or\ifcase#2\or\or\or\or\or\or
\or\qbezier(-1,0)(-0.516,-0.138)(-0.866,-0.5)
\or\qbezier(-1,0)(-0.45,-0.26)(-0.5,-0.866)
\or\qbezier(-1,0)(-0.327,-0.327)(0,-1)
\or\qbezier(-1,0)(-0.179,-0.311)(0.5,-0.866)
\or\qbezier(-1,0)(-0.0536,-0.2)(0.866,-0.5)
\fi
\or\ifcase#2\or\or\or\or\or\or\or
\or\qbezier(-0.866,-0.5)(-0.378,-0.378)(-0.5,-0.866)
\or\qbezier(-0.866,-0.5)(-0.26,-0.45)(0,-1)
\or\qbezier(-0.866,-0.5)(-0.12,-0.446)(0.5,-0.866)
\or\qbezier(-0.866,-0.5)(0,-0.359)(0.866,-0.5)
\fi
\or\ifcase#2\or\or\or\or\or\or\or\or
\or\qbezier(-0.5,-0.866)(-0.138,-0.516)(0,-1)
\or\qbezier(-0.5,-0.866)(0,-0.52)(0.5,-0.866)
\or\qbezier(-0.5,-0.866)(0.12,-0.446)(0.866,-0.5)
\fi
\or\ifcase#2\or\or\or\or\or\or\or\or\or
\or\qbezier(0,-1.)(0.138,-0.516)(0.5,-0.866)
\or\qbezier(0,-1.)(0.26,-0.45)(0.866,-0.5)
\fi
\or\ifcase#2\or\or\or\or\or\or\or\or\or\or
\or\qbezier(0.5,-0.866)(0.378,-0.378)(0.866,-0.5)
\fi\fi\fi\fi}
\def\FullChord[#1,#2]{
\Endpoint[#1]
\Endpoint[#2]
\Arc[#1]
\Arc[#2]
\Chord[#1,#2]
}
\def\EndChord[#1,#2]{
\Endpoint[#1]
\Endpoint[#2]
\Chord[#1,#2]
}
\def\Picture#1{
\begin{picture}(2,1)(-1,-0.167)
#1
\end{picture}
}
\def\DottedChordDiagram[#1,#2]{
\Picture{\DottedCircle \FullChord[#1,#2]}
}
\newcommand{\e}{\varepsilon} 
\newtheorem{theor}{Theorem}[section] 
\newtheorem{thm}[theor]{Theorem} 
 
\newtheorem{prop}[theor]{Proposition} 
\newtheorem{corol}[theor]{Corollary} 
\newcommand{\begproof}{\noindent \\ {\em Proof\/.} } 
\newcommand{\eproof}{ {\quad \ \hfill $\Box$ \\ \medskip}} 
\renewcommand{\a}{\alpha} 
\renewcommand{\b}{\beta} 
\newcommand{\g}{{\gamma}} 
\renewcommand{\d}{\delta}
\newcommand{\s}{\sigma} 
\newcommand{\be}{\begin{equation}} 
\newcommand{\ee}{\end{equation}} 
\newcommand{\bea}{\begin{eqnarray}} 
\newcommand{\eea}{\end{eqnarray}} 
\newcommand{\bean}{\begin{eqnarray*}} 
\newcommand{\eean}{\end{eqnarray*}} 
\newcommand{\text}[1]{ {\rm {#1} } }           
\newcommand{\Bbb}{\bf}                       
\newcommand{\mathbf}{\bf}                    
\newcommand{\mathrm}{\rm}                    
\newcommand{\Q}{{\Bbb Q}} 
\newcommand{\C}{{\Bbb C}} 
\newcommand{\R}{{\Bbb R}}

\newcommand{\D}{\Delta}

\newcommand{\cF}{{\cal F}} 
\newcommand{\cK}{{\cal K}} 
\newcommand{\cD}{{\cal D}}

\newcommand{\cU}{{\cal U}} 
\newcommand{\cV}{{\cal V}} 
\newcommand{\cW}{{\cal W}} 
\newcommand{\Z}{{\Bbb Z}} 

\newcommand{\BN}{Bar-Natan} 
\newcommand{\Gar}{Garoufalidis}
\newcommand{\AC}{Alexander}
\newcommand{\Fn}{Feynman} 
\newcommand{\Fd}{\Fn\ diagram} 
\newcommand{\gl}{\mbox{$gl(1|1)$}} 
\newcommand{\hW}{\widehat{W}}
\newcommand{\ra}{\rightarrow} 
 
\newcommand{\Vas}{Vassiliev} 
\newcommand{\VI}{\Vas\ invariants} 
\newcommand{\ki}{knot invariant} 
\newcommand{\ws}{weight system} 

\newcommand{\half}{\frac{1}{2}}
\newcommand{\ov}{\overline}

\newcommand{\MM}{Melvin-Morton}

\begin{document}

\title{
\bf Universal \ws s and the \MM\ expansion of the
colored Jones knot invariant}

\author{
Arkady Vaintrob
\\ {\em Department of Mathematics, University of Utah%
\thanks{on leave from New Mexico State University, Las
Cruces, NM 88003}%
,} 
\\ {\em Salt Lake City, UT 84112, USA \qquad e-mail: {\tt
vaintrob@math.utah.edu} }
}

\date{}
\maketitle

 
\begin{abstract}  
We study the asymptotic expansion  of the colored Jones
polynomial (the \MM\ expansion)
using  a recursive formula for the deframed universal
\ws\ for the $sl_2$ Lie algebra. 
Combined with the formula
for the universal \ws\ for the Lie superalgebra \gl\ (which
corresponds to the \AC\ knot polynomial) this formula gives a
very short proof of the \MM\ conjecture relating the colored
Jones invariant and the \AC\  polynomial of knots.  
\end{abstract} 
 

\section{Introduction}

The two most famous invariants of oriented links in $S^3$, the Alexander
and Jones polynomials, have very similar combinatorial definitions.

Let  $U$ be the unknot and $K_+$, $K_-$, and $K_=$ three knots (or
links) identical outside of a small ball whose intersection with this ball
looks as shown on Fig.~\ref{fig:skein}.

\begin{figure}[hbt]
\def\KZero{
\begin{picture}(2,2)(-1,-1)
\put(0,0){\circle{2}}
\put(-0.707,-0.707){\vector(1,1){1.414}}
\put(0.707,-0.707){\vector(-1,1){1.414}}
\put(0,0){\circle*{0.15}}
\end{picture}}

\def\KPlus{
\begin{picture}(2,2)(-1,-1)
\put(0,0){\circle{2}}
\put(-0.707,-0.707){\vector(1,1){1.414}}
\put(0.707,-0.707){\line(-1,1){0.6}}
\put(-0.107,0.107){\vector(-1,1){0.6}}
\end{picture}}

\def\KMinus{
\begin{picture}(2,2)(-1,-1)
\put(0,0){\circle{2}}
\put(-0.707,-0.707){\line(1,1){0.6}}
\put(0.107,0.107){\vector(1,1){0.6}}
\put(0.707,-0.707){\vector(-1,1){1.414}}
\end{picture}}

\def\KII{
\begin{picture}(2,2)(-1,-1)
\put(0,0){\circle{2}}
\qbezier(-0.707,-0.707)(0,0)(-0.707,0.707)
\qbezier(0.707,-0.707)(0,0)(0.707,0.707)
\put(-0.607,0.607){\vector(-1,1){0.1414}}
\put(0.607,0.607){\vector(1,1){0.1414}}
\end{picture}}

\begin{displaymath}
\mathop{\KPlus}_{K_+}~,\qquad
\mathop{\KMinus}_{K_-}~,\qquad
\mathop{\KII}_{K_=}~,\qquad
\mathop{\KZero}_{K_0}~.
\end{displaymath}
\caption{}
\label{fig:skein}
\end{figure}

The Alexander polynomial $\D$ is the link invariant uniquely
determined by the following conditions (skein relations)
\be
\D(K_+,t) - \D(K_-,t) =(t^{\half} - t^{-\half})\D(K_=,t), \qquad \D(U,t)=1~,
\label{eq:skein1}
\ee
and the Jones polynomial of links is defined by the relations
\be
t^{-1}V(K_+,t) - tV(K_-,t) =(t^{\half} - t^{-\half})V(K_=,t)~, \qquad
V(U,t)=1~. 
\label{eq:skein2}
\ee

Both $\D$ and $V$ are Laurent polynomials in $\sqrt{t}$,
but if $K$ has only one component (i.e. $K$ is a knot),
$\D(K)$ and $V(K)$ belong to $\C[t,t^{-1}]$.

Despite the similarity between the skein relations for the Jones and
Alexander invariants,  
there is a huge gap in our understanding of their nature. 
The Alexander polynomial has a solid   topo\-logical meaning (as a
determinantal invariant  of the $\Z[t,t^{-1}]$-module
$H^1(\widetilde{M})$, where 
\mbox{ $\widetilde{M} \to S^3 \setminus K$ }
is the infinite cyclic cover of the knot complement corresponding to
the commutator  subgroup $G'=[G,G]$ of the knot group $G=\pi_1(S^3
\setminus K)$). 

The Jones polynomial, on the other hand, is defined purely
combinatorially and even though it was used to settle some very old
conjectures in knot theory, its topological meaning is still very
obscure.

The theory of quantum groups offered a new viewpoint on 
the Jones polynomial. It was found that 
$V(K)$ can be  constructed by working with the standard
two-dimensional representation of the quantum group $sl_{2,q}$. 
An arbitrary, $d$-dimensional, representation of $sl_{2,q}$
gives rise to a generalization of $V$, 
the so-called {\em colored Jones invariant\/} $V^d$
(see \cite{KR,RT} for more details).

Morton and Strickland \cite{MS} found 
a combinatorial formula for computing $V^d(K,t)$
proving that $V^d$ is determined by all
the cablings of the Jones polynomial.

In \cite{MM} Melvin and Morton  studied pieces of $V^d(K,t)$ appearing
in a certain power series expansions and conjectured
that the Alexander polynomial of a knot is 
determined by its colored Jones invariant.
They proved that the coefficients $v_{in}$
of the power series expansion 
\be
V^d(K,e^z) = \sum_{i,n \ge 0} v_{in}(K)d^iz^n
\label{eq:expansion}
\ee
of $V^d$ in variables 
$d$ and $z=\log t$
are \Vas\ \ki s
of order $\le n$ and that $v_{in}=0$ if $i>2n$.
They formulated the following {\em \MM\ conjecture\/}
(later Morton~\cite{M} proved it for torus knots).  

\begin{thm} \label{thm:mm} \ \nopagebreak

(i) The coefficient matrix $(v_{in})$ of the expansion
(\ref{eq:expansion}) is ``lower triangular,'' i.e.
\be
 v_{in} = 0 \quad \text{for} \ i>n 
\label{eq:MM0}
\ee

(ii)  the leading (``diagonal'') term 
\be
V_0(K,z) = \sum_{n\ge 0} v_{nn}(K) z^n 
\label{eq:V0}
\ee
is the inverse of the renormalized Alexander polynomial of $K$
\be
V_0(K,z) \cdot \frac{z}{e^{z/2} - e^{-z/2}} \Delta(K,e^z) = 1~.
\label{eq:MM}
\ee
\end{thm}

Rozansky \cite{Roz} derived relations~(\ref{eq:MM0}) and~(\ref{eq:MM})
from Witten's~\cite{W} path integral interpretation of
the Jones invariant  using quantum field theory tools.

Finally, \BN\ and  \Gar\ \cite{BNG} proved the \MM\ conjecture by
finding nice combinatorial expressions for the \ws s corresponding to
the leading terms $v_{nn}$ of the expansion~(\ref{eq:expansion}) and the 
coefficients of the \AC\ polynomial.  
Their method cannot be used to study
off-diagonal coefficients of~(\ref{eq:expansion}).

One of the difficulties in studying the colored Jones invariant is
that it deals with all representations of $sl_2$ at once. In this
paper we suggest a different approach to analyzing
the \MM\ expansion~(\ref{eq:expansion}) for which this problem does
not arise. Our method is based on universal \VI\ --- \Vas\ \ki s
corresponding to \ws s with values in the center $Z(\cU(L))$
of the universal enveloping algebra of a Lie (super)algebra $L$. 
The $sl_2$ universal \Vas\ invariant is a power series with
coefficients in $Z(\cU(sl_2))= \C[c]$, where
$c$ is the quadratic Casimir. It is equivalent to the colored Jones
invariant, since in the $d$-dimensional representation  $c$ acts as
the scalar $(d^2-1)/2$. Thus, we can study the colored Jones invariant 
entirely in terms of $sl_2$ alone without working with its
different representations. 

Very little is known about universal \VI\ and the
corresponding \ws s in general. The only non-trivial examples with
complete answers are the cases of the Lie algebra $sl_2$ and of the
Lie superalgebra $gl(1|1)$ which correspond respectively to the 
colored Jones and the Alexander \ki s.

Chmutov and Varchenko \cite{CV} found a recursion formula for the 
universal $sl_2$ invariant for {\em framed\/} knots. In this paper
we derive a similar relation for the $sl_2$ invariant of 
{\em unframed\/} knots (i.e. of the colored Jones invariant)
that enables us to compute recursively the \ws s corresponding to all
the ``lines'' 
$$
   V_p=\sum\limits_{n\ge 0} v_{n,n+p} z^n 
$$ 
of the expansion~(\ref{eq:expansion}).  

The Alexander polynomial is, essentially, the
universal invariant for the simplest interesting Lie superalgebra
$\gl$ (cf.~\cite{KS,Va,VaiYB,FKV}). In~\cite{FKV} we studied the
universal $\gl$ invariants (both for framed and unframed knots) and
found recursion formulas for the corresponding \ws s.

The recursive relation for the deframed $sl_2$-invariant looks very
similar to the one for $\gl$. This gives a very simple proof of
Theorem~\ref{thm:mm} and paves a way to a better understanding the 
lower lines $V_p$ in~(\ref{eq:expansion}).

The proofs of the recursion formulas for both the $sl_2$ and $\gl$
universal invariants are based on two fundamental relations 
(\ref{eq:lagrange}) and (\ref{eq:gl11fund}) between
invariant tensors of order four on these Lie (super)algebras.

The relation for $sl_2$ has been known for a long
time. In the $so_3$ realization of $sl_2$ it is the famous Lagrange's
identity of vector algebra in $\R^3$:
$$\mathbf 
[a\times b]\cdot [c\times d] = (a\cdot c)(b\cdot d) - (a\cdot
d)(b\cdot c)~.
$$
The corresponding relation for  $gl(1|1)$ was discovered in~\cite{FKV}.

These relations, in fact, are responsible for the skein
relations~(\ref{eq:skein1}) and~(\ref{eq:skein2}) and for  
other properties of the Jones and Alexander \ki s.

The plan of the paper is as follows. In Section~2 we collect necessary
information on \Vas\ \ki s and their relations with Lie algebra-type
structures. In particular, we explain, following~\cite{BNG}, why it is
enough to establish~(\ref{eq:MM0}) and~(\ref{eq:MM}) on the level of
\ws s. In Section~3 we discuss the \ws s for  the coefficients of the
\MM\ expansion~(\ref{eq:expansion}) of the colored Jones invariant.
In particular, we prove a recursion formula for the deframed universal
$sl_2$ \ws. In Section~4 we review the results of~\cite{FKV} on the
universal $\gl$ \Vas\ invariant and its relation with the Alexander
polynomial. Then we show that Theorem~\ref{thm:mm} is a simple
corollary of our formulas for the universal $sl_2$ and $\gl$
invariants.

\section{\VI\ and universal \ws s}

Here we review some facts about \VI\
and their relationship with Lie algebra-type structures. For more
details see~\cite{BN,Konts,VaiYB}.
 
\subsection{Vassiliev invariants, chord diagrams and weight systems}

A {\em singular knot\/} is an immersion $K: S^1 \rightarrow \R^3$ with
a finite number of double self-intersections with distinct tangents.
The set of singular knots with $n$ double points is denoted
by $\cK_n$.

A {\em chord diagram\/} of order $n$ is an oriented circle with $n$
disjoint pairs of points ({\em chords\/}) on it up to an orientation
preserving diffeomorphism of the circle.  Denote by $\cD_n$ the set of
all chord diagrams with $n$ chords.

Every singular knot $K\in \cK_n$ has a chord diagram $ch(K)\in \cD_n$
whose chords are the inverse images of the double points of $K$.

Every  knot invariant $I$ with values in an abelian
group $k$ extends to an invariant of singular knots by the rule
\be 
I(K_0) = I(K_+) - I(K_-),    \label{eq:vasrel}
\ee
where $K_0$, $K_+$, and $K_-$ are singular knots which differ only
inside a small ball as shown on Fig.~\ref{fig:skein}.

A knot invariant $I$ is called an {\em invariant of order\/} ($\le$)
$n$ if $I(K)=0$ for any $K\in \cK_{n+1}$.

The $k$-module of all invariants of order $n$ is denoted by $\cV_n$.  We
have the obvious filtration
\begin{displaymath}
\cV_0 \subset \cV_1 \subset \cV_2 \ldots
\subset \cV_n \subset \ldots~.
\end{displaymath}

Elements of $\cV_\infty = \bigcup\limits_n \cV_n$ 
are called {\em invariants of finite type\/} or {\em \VI}. 

Similarly the space $\cV_n^f$ of \VI\ of framed knots  is defined. 

\medskip

An immediate corollary of the definition of Vassiliev invariants is
that the value of an invariant $I \in V_n$ on a singular knot $K$ with
$n$ self-intersections depends only on the diagram $ch(K)$ of $K$.  In
other words, $I$ descends to a function on $\cD_n$ which we still
denote by $I$.
These functions satisfy two groups of relations
\begin{equation}
I\left(\Picture{\DottedCircle\FullChord[2,10]\Arc[1]\Arc[0]\Arc[11]}\right)
= 0\label{eq:1term}
\end{equation}
and
\begin{equation}
I\left(
\Picture{\DottedCircle\FullChord[1,8]\Arc[2]\FullChord[5,9]}
\right) -
I\left(
\Picture{\DottedCircle\FullChord[1,9]\Arc[2]\FullChord[5,8]}
\right) +
I\left(
\Picture{\DottedCircle\FullChord[2,5]\Arc[8]\FullChord[1,9]}
\right) -
I\left(
\Picture{\DottedCircle\FullChord[1,5]\Arc[8]\FullChord[2,9]}
\right) = 0.\label{eq:4term}
\end{equation}
\ 

\medskip

Besides the explicitly shown chords, the diagrams may contain other
chords with endpoints lying on the dotted arcs, provided that these
chords  are the same in all diagrams appearing in the same relation.
\medskip

A function $W: {\cal D}_n \rightarrow k$ is called a {\em weight
system of order $n$\/} if it satisfies the
{\em four-term relations\/}~(\ref{eq:4term}). If, in addition, $W$
satisfies the {\em one-term relations\/}~(\ref{eq:1term}),  we
call it a {\em strong \ws}.   

Denote by $\cW_n$ (resp.\ by $\ov{\cW}_n$) the set of all weight
systems (resp.\ strong weight systems) of order $n$.

A \Vas\ invarant (resp.\ a \Vas\ invariant of framed knots) of order
$n$ defines a strong \ws\ (resp.\ a \ws) of order $n$ and 
it is easy to see that the natural maps
 $\cV^f_n/\cV^f_{n-1} \rightarrow \cW_n$ and
 $\cV_n/\cV_{n-1} \rightarrow \ov{\cW}_n$ 
are injective.

 The remarkable fact proved by Kontsevich~\cite{Konts} is that these
 maps are also surjective (at least when $k \supset \Q$).  In other
 words, each (strong) \ws\ of order $n$ is a restriction to
 ${\cal{D}}_n$ of some Vassiliev invariant.
 
Kontsevich proved that 
 $$
\cV_n/\cV_{n-1} \simeq \ov{\cW}_n~,
   \ \text{and} \ 
\cV^f_n/\cV^f_{n-1} \simeq \cW_n~
$$
by explicitly constructing 
splitting maps
 \be
Z: \cW_n \to \cV^f_n \quad \text{and} \quad \bar{Z}: \ov{\cW}_n \to
\cV_n
\label{eq:splitting}
\ee

If $k$ is a commutative ring, then the
product of two \VI\ of orders $m$ and $n$ is a \Vas\ invariant 
of order $m+n$, therefore $\cV$ is a filtered algebra.

The space  
$$
\bigoplus\limits_{n \ge 0}\cW_n =
\bigoplus\limits_{n \ge 0}\cV^f_{n+1}/\cV^f_n
$$
becomes the adjoint graded algebra of $\cV$
with the product defined as follows.
 For
$W\in \cW_m$, $W'\in \cW_n$ and $D\in \cD_{m+n}$ we have 
\be
(W\cdot W')(D) = \sum_{{E\subset D,} \atop {|E|=m}}
W(E)W'(D\setminus E).
\label{eq:wsproduct}
\ee

Many \ki s such as the Alexander and Jones polynomials are not \VI,
whereas their coefficients (after an appropriate change of variables)
are. Therefore, we can associate with such invariants not just one
\ws, but a sequence of \ws s $w_0, w_1, \ldots, w_n, \ldots$
where $w_n \in \cW_n$. 
We call elements of $\cW = \prod\limits_{n \ge 0} \cW_n$ {\em \Vas\
series\/} and write them as formal sums 
 $W=w_0+w_1+w_2 + \ldots$ (or sometimes as 
formal power series $\sum w_nz^n$, where $z$ is a formal parameter).
\Vas\ series can be viewed as linear functionals 
on the space of diagrams
$\cD=\bigoplus\limits_{n \ge 0} \cD_n$:
for $D\in \cD_n$ and $W \in \cW$ we define $W(D) =w_n(D)$. 
The space $\cW$ also has an algebra structure
\be
(W\cdot W') (D) = \sum_{E\subset D} W(E)W'(D\setminus E).
\label{eq:product}
\ee

Every \ws\ of order $n$ gives a strong \ws\ of order $n$ by 
means of a canonical projection ({\em deframing\/})
\be
p: \cW_n \to \ov{\cW}_n
\label{eq:defr}
\ee
given by
\be
p(W)(D) = \sum_{E\subset D} (-1)^{|E|}
W\Bigl(\Theta^{|E|}\cdot(D\setminus E)\Bigr),
\label{eq:defr1}
\ee
where $\Theta$ is the unique chord diagram with one chord
and the product of two chord diagrams $D_1$ and $D_2$ is just their
connected sum $D_1 \cdot D_2$. The product of diagrams is well-defined modulo
 four-term relations~(\ref{eq:4term}). 

A \Vas\ series $W \in \cW$ is called {\em multiplicative\/} if
$$ 
  W(D_1\cdot D_2) = W(D_1)\cdot W(D_2) \ \text{for\ any}\ D_1, D_2 \in \cD.
$$
Combining equations~(\ref{eq:product}) and~(\ref{eq:defr1})
we obtain a  convenient formula for deframing multiplicative \ws s.

\begin{prop} 
Let $W$ be a multiplicative \Vas\ series and
 $U$ be the \Vas\ series 
$$
U(D) = (-c)^{|D|}, \  \text{where} \  c=W(\Theta).
$$
Then the deframed strong \ws\ $\ov{W} = p(W)$
is given by 
\be
\ov{W} = W\cdot U \quad \text{\em or} \quad
\ov{W}(D) = \sum_{E\subset D} (-c)^{|E|}\cdot W(D\setminus E).
\label{eq:deframed}
\ee
\end{prop}
\eproof

\subsection{Weight systems coming from Lie algebras} 
\label{sec:FG}

Here we recall a construction that
assigns a family of weight systems to every Lie
(super)algebra with an invariant inner product. 

First, we introduce a more general class of diagrams.
\medskip

A {\em \Fd \/} of order $p$ is a graph with $2p$ vertices of degrees 1
or 3 with a cyclic ordering on the set of its univalent ({\em
external\/}) vertices and on each set of three edges meeting at a
trivalent ({\em internal\/}) 
vertex.\footnote{\Fd s
are called  Chinese
Character diagrams in~\cite{BN}, but they are indeed \Fd s
arizing in the perturbative Chern-Simons-Witten
quantum field theory~\cite{W,BN2}.}   
Let $\cF_p$ denote the set of all
\Fd s with $2p$ vertices (up to the natural equivalence of graphs with
orientations). The set $\cD_p$ of chord diagrams with $p$ chords is a
subset of $\cF_p$.

\medskip

Any \ws\ $W:\cD_p \to k$ can be extended to $\cF_p$ by the following
rule, the {\em three-term relations\/}

\begin{equation}\label{eq:3term}
W\Picture{
\Arc[9]
\Arc[10]
\Arc[11]
\Arc[0]
\DottedArc[1]
\DottedArc[8]
\thinlines
\put(0.7,-0.7){\circle*{0.15}}
\qbezier(-0.1,0.1)(0.2,-0.5)(0.5,0.1)
\qbezier(0.2,-0.2)(0.3,-0.6)(0.7,-0.7)
\put(0,-1.4){\makebox(0,0){${}_{D_Y}$}}
}\quad = \quad
W\Picture{
\Arc[9]
\Arc[10]
\Arc[11]
\Arc[0]
\DottedArc[8]
\DottedArc[1]
\Endpoint[10]
\Endpoint[11]
\thinlines
\put(0.5,-0.866){\line(-3,5){0.58}}
\put(0.866,-0.5){\line(-3,5){0.36}}
\put(0,-1.4){\makebox(0,0){${}_{D_{\mid\mid}}$}}
} \quad - \quad 
W \ \Picture{
\Arc[9]
\Arc[10]
\Arc[11]
\Arc[0]
\DottedArc[8]
\DottedArc[1]
\Endpoint[10]
\Endpoint[11]
\thinlines
\put(0.5,0.1){\line(0,-1){0.97}}
\put(0.866,-0.5){\line(-5,3){1}}
\put(0,-1.4){\makebox(0,0){${}_{D_X}$}}
}~.
\end{equation}

\ 

\bigskip

The space
 $\langle\cF_p\rangle /\langle D_Y-D_{\mid\mid}+D_X\rangle$ 
of \Fn\ diagrams modulo three-term relations is isomorphic 
to the space generated by $\cD_p$ modulo four-term relations
(\ref{eq:4term}), i.e.\ it is canonically dual to $\cW_p$.

In addition,
the following local relations hold for internal vertices of \Fd s
\be
W 
\begin{picture}(2,2)(0,0.375)
\qbezier(0.5,2)(1.65,1.3)(1.5,1)
\qbezier(1.5,2)(0.35,1.3)(0.5,1)
\qbezier(0.5,1)(1,0)(1.5,1)
\put(1,0){\line(0,1){0.5}}
\end{picture}
\quad = - 
W
\begin{picture}(2,2)(0,0.375)
\qbezier(0.5,2)(1,0)(1.5,2)
\put(1,0){\line(0,1){1}}
\end{picture}
\ {\mathrm and}\quad 
W
\begin{picture}(2,2)(-1,-1)
\put(0,-1.4){\line(1,1){1.4}}
\qbezier(-0.1,0.1)(0.2,-0.5)(0.5,0.1)
\qbezier(0.2,-0.2)(0.3,-0.6)(0.7,-0.7)
\end{picture} \quad =
W
\begin{picture}(2,2)(-1,-1)
\put(-0.034,-1.4){\line(1,1){1.4}}
\put(0.5,-0.866){\line(-3,5){0.58}}
\put(0.866,-0.5){\line(-3,5){0.36}}
\end{picture} \quad - 
W
\begin{picture}(2,2)(-1,-1)
\put(-0.034,-1.4){\line(1,1){1.4}}
\put(0.5,0.1){\line(0,-1){0.97}}
\put(0.866,-0.5){\line(-5,3){1}}
\end{picture}
.
\label{eq:asjac}
\ee

\bigskip

Let $L$ be a Lie (super)algebra with an $L$-invariant inner product
$b: L\otimes L \ra k$. To each \Fn\ diagram $F$ with $m$ univalent
vertices  we assign a tensor 
$$
   T_{L,b}(F) \in L^{\otimes m}
$$ 
as follows.

The Lie bracket $[\ ,\ ]: L\otimes L \to L$ can be considered
as a tensor in $L^*\otimes L^* \otimes L$. The inner product 
$b$ allows us to identify the $L$-modules 
$L$ and $L^*$, and therefore  $[\ ,\ ]$ can be considered as a tensor
$f\in (L^*)^{\otimes 3}$ and $b$ gives rise to an invariant symmetric
tensor $c \in L\otimes L$.

For a \Fn\ diagram $F$ denote by $T$ the set of its trivalent
vertices, by $U$ the set of its univalent (exterior) vertices, and by
$E$ the set of its edges. Taking $|T|$ copies of the tensor $f$ and $|E|$
copies of the tensor $c$ we consider a new tensor
$$
\widetilde{T}_L(F) = \Bigl(\bigotimes_{v\in T} f_v\Bigr) \otimes
                     \Bigl(\bigotimes_{\ell\in E} c_\ell\Bigr)
$$
which is an element of the tensor product
$$
{\cal{L}}^F = 
\Bigl( \bigotimes_{v \in T}
(L^*_{v,1}\otimes L^*_{v,2} \otimes L^*_{v,3})
\Bigr)
\otimes
\Bigl( \bigotimes_{\ell \in E}
(L_{\ell,1}\otimes L_{\ell,2})
\Bigr),
$$
where $(v,i), \ i=1,2,3,$ mark the three edges meeting at the
vertex $v$ (consistently with the cyclic ordering of these edges), and
$(\ell,j), \ j=1,2,$ denote the endpoints of the edge $\ell$. 

Since $c$ is symmetric and $f$ is completely antisymmetric,
the tensor $\widetilde{T}_L(F)$ does not depend on the choices of
orderings.

If $(v,i)=\ell$ and $(\ell,j)=v$, there is a natural contraction map
$L^*_{v,i} \otimes L_{\ell,j} \to k$. Composition of all such 
contractions gives us a map
$$
{\cal{L}}^F \longrightarrow  \bigotimes_{u \in U} L = L^{\otimes m}, \ \text{where} \ m=|U|.
$$

The image of $\widetilde{T}_L(F)$ in
$L^{\otimes m}$ is denoted by $T_{L,b}(F)$ (or usually just by
$T_{L}(F)$).

(In the case of Lie superalgebras we also have to take special care of
signs.  See \cite{VaiYB} for details.)

For example, for the diagrams 
\def\DiagC{
\begin{picture}(3,1)
\qbezier(0.5,0)(1.5,2)(2.5,0)
\thicklines
\put(0,0){\line(1,0){3}}
\end{picture}
}
\def\Bubble{
\begin{picture}(4,1.5)
\put(1.5,0){\oval(2,2)[tl]}
\put(2.5,0){\oval(2,2)[tr]}
\put(2,1){\circle{1}}
\thicklines
\put(0,0){\line(1,0){4}}
\end{picture}}
\def\DiagramK{
\begin{picture}(5,1.5)
\qbezier(0.5,0)(1.25,1.5)(2,0)
\qbezier(3,0)(3.75,1.5)(4.5,0)
\qbezier(1.25,0.75)(2.5,2)(3.75,0.75)
\thicklines
\put(0,0){\line(1,0){5}}
\end{picture}}
\be
C= \DiagC~, 
\quad 
B = \Bubble~, 
\quad \text{and} 
\quad
K =\DiagramK
\label{eq:DiagK}
\ee
we have 
$$T_L(C) = \sum_{ij}  b^{ij}e_i\otimes e_j= c,
$$ 
the Casimir element corresponding to the inner product~$b$, \
$$
T_L(B) = \sum b^{is}b^{tj}b^{kp}b^{lq}f_{skl}f_{pqt}e_i \otimes e_j,
$$
the tensor in $L\otimes L$
corresponding to the Killing form on $L$ under the identification 
$L^* \simeq L$, and
$$
T_L(K) = \sum b^{in} b^{jp} b^{qr} b^{kt} b^{\ell s}
f_{npq} f_{tsr}^k e_i \otimes e_j
\otimes e_k \otimes e_\ell~,
$$
 where $f^i_{jk}$ are the structure constants of 
$L$ in a basis $e_1,e_2,\ldots$~.

Tensor $T_L(F)$ is invariant with respect to the $L$-action on
$L^{\otimes m}$ and its image $W_L(F)$ in the universal enveloping
algebra $\cU(L)$ belongs to the center 
$Z(\cU(L))=\cU(L)^L$ and does not depend on the place where we
cut the Wilson line to obtain a linear ordering of the external
vertices of $F$.   

The conditions (\ref{eq:asjac}) and (\ref{eq:3term}) are automatically
satisfied for $W_L$: the relations (\ref{eq:asjac}) are 
the anticommutativity 
and the Jacobi identity for the Lie bracket, and (\ref{eq:3term}) is
just the definition of the universal enveloping algebra as a quotient
of the tensor algebra of $L$.

Therefore, for every Lie algebra $L$ with an invariant inner product
there exists a natural \Vas\ series
 $W_{L}: \cD \rightarrow Z(\cU(L))$.

The \Vas\ series 
$$
W_L: \cD \to Z(\cU(L)) 
$$
is called {\em the universal weight system}
corresponding to $L$ (we are suppressing  $b$ from the notation).
It is universal in the sense that any \Vas\  series $W_{L,R}$
constructed using a representation $R$ of the Lie algebra $L$ (see
\cite{BN})  is an evaluation of $W_L$:
$$
W_{L,R}(D) = \text{Tr}_{R}\bigl(W_L(D)\bigr).
$$

By its construction, the universal
 \Vas\ series $W_L$ is multiplicative, i.e. 
$W_L(D_1 \cdot D_2) = W_L(D_1)W_L(D_2)$.

\subsection{Reduction of the \MM\ conjecture to \ws s}
\label{sec:canon}

Given a semi-simple Lie algebra $L$  with an invariant inner 
product $b$ and a representation $R$, 
there are two ways to construct knot invariants from this data.
First, we can use the Reshetikhin-Turaev construction~\cite{RT} 
based on quantum groups to get invariants $I_{L,R}$ and $\bar{I}_{L,R}$
(of framed and unframed knots, resp.)\ with values
in $\Z[t,t^{-1}]$. Second, we can apply Kontsevich's splitting map
to \Vas\ series $W_{L,R}=\text{Tr}_R(W_L)$ and $\ov{W}_{L,R}$ discussed 
in the previous section. 

The remarkable fact is that these two 
constructions are equivalent. As it was first noticed by Birman and 
Lin~\cite{BL}, the coefficients of the power series expansions of 
$\bar{J}_{L,R}(z)=\bar{I}_{L,R}(e^z)$ are \VI. 
Piunikhin~\cite{P} proved that the corresponding series of
\ws s coincide with $\ov{W}_{L,R}$. Kassel~\cite{Kassel} and
Le and Murakami~\cite{LM} showed that, conversely, Kontsevich's 
construction applied to $\ov{W}_{L,R}$ (resp.\ to $W_{L,R}$)
gives the sequence of the coefficients of $\bar{J}_{L,R}$
(resp.\ $J_{L,R}$). 

A \Vas\ invariant is called {\em canonical\/} 
if it belongs to the image of 
Kontsevich's map $Z: \cW_n \to \cV_n$~.
A formal power series  
$\sum\limits_{n \ge 0} v_n z^n \in \cV[[z]]$ 
is called {\em canonical\/} if every coefficient $v_n$ 
is a canonical \Vas\ invariant of order $\le n$.

The colored Jones invariant is the canonical invariant
$Z(\frac{1}{d}\ov{W}_{sl_2,R_d})$ corresponding to the $d$-dimensional
representation of $sl_2$. Since a canonical invariant is uniquely
determined by its \ws, to prove part (i) of the \MM\ 
conjecture~(\ref{thm:mm}) it is enough to show that the
weight system of $v_{in}$ in (\ref{eq:expansion})
vanishes for $i > n$.

\BN\ and \Gar~\cite{BNG} proved that
\be
\tilde{\D}(z) =  \frac{z}{e^{z/2} - e^{-z/2}}\D(e^z),
\label{eq:alexnorm}
\ee
where $\D(K,t)$ is the Alexander polynomial of knot $K$,
is a canonical series. (This also follows from the
generalization of the Kassel-Le-Murakami theorem for classical
Lie superalgebras and from the 
fact that the \AC\ invariant and the corresponding \Vas\
series $C$ come from \gl, see \cite{KS,Va,FKV}). 

The product of two canonical invariants or \Vas\ series
$Z(W_1)$ and $Z(W_2)$ is again canonical with the \ws\
equal to $W_1W_2$.
Therefore, to prove part (ii) of the \MM\ conjecture
it is enough to establish that the \Vas\ series corresponding to
the two factors in (\ref{eq:MM}) are inverses of each other.

Theorem~\ref{thm:mm} is now reduced to the following
relations on the level of \ws s.

\begin{prop} \label{prop:ws}
Let 
\be
{\hW}_{sl}=\frac{1}{d}\sum_{n \ge 0} \ov{W}_{sl_2,R_d,n}
=   \sum_{i,n \ge 0} w_{in} d^i
\label{eq:coefmm}
\ee
be  the deframed \Vas\ series (\ws) coming
from the $d$-dimensional representation of $sl_2$ with the standard
metric $\langle x,y \rangle = \text{Tr}(xy)$
(normalized by dividing by $d$).
Then
\be (i) \qquad w_{in}=0 \  \text{~for~} \  i>n
\label{eq:mm1}
\ee
and
\be (ii) \qquad (\sum_{n \ge 0} w_{nn} )\cdot C = \e,
\label{eq:mm2}
\ee
where $C = \sum C_n z^n \in \cW$ is the \Vas\ series
of the normalized \AC\ invariant $\tilde{\D}$ (\ref{eq:alexnorm})
and $\e$ is the \Vas\ series 
$$
\e(D) = \cases {1 & if $|D|=0$ ,\cr 0 & if $|D| > 0$~. }
$$
\end{prop}
\medskip

\BN\ and \Gar~\cite{BNG} proved~(\ref{eq:mm2}) by
finding combinatorial  formulas for both \ws s in the
left-hand side of this equation. We will show that
both~(\ref{eq:mm1}) and~(\ref{eq:mm2}) are simple corollaries
of the recursion formulas for the $sl_2$ and $\gl$
universal \ws s.

\section{Universal $sl_2$ \ws s}

The standard choice of an invariant metric on $L=sl_2$ is 
$$\langle x,y \rangle = \text{Tr}_V (xy). $$
This metric is equal to one fourth of the  Killing form 
$\langle x,y \rangle_K = \text{Tr}(ad_x ad_y),$
 therefore 
\be
W_L\ \Bubble = 4 W_L \DiagC~.
\label{eq:slkilling}
\ee

The center of $\cU(L)$ is $\C[c]$ where $c$ is the quadratic Casimir --- the element of
$S^2(L)$ corresponding to the metric under the
identification $L^* \simeq L$.

In this setting 
$$
W(D) = c^n - 2p c^{n-1} + \text{terms\ of\ lower\ order\ in\ } c,
$$
where $n=|D|$ is the number of chords in the diagram $D$, and $p$ is
the number of pairs of intersecting chords in $D$.

The following relation between invariant tensors in $L^{\otimes 4}$ is responsible
for all properties of the \Vas\ series $W_L$.

Let 
\begin{displaymath}
P = \begin{picture}(4,1)
\qbezier(0.5,0)(1.5,2)(2.5,0)
\qbezier(1.5,0)(2.5,2)(3.5,0)
\thicklines
\put(0,0){\line(1,0){4}}
\end{picture} - 
\begin{picture}(4,1.5)
\qbezier(0.5,0)(2,2.5)(3.5,0)
\qbezier(1.5,0)(2,1)(2.5,0)
\thicklines
\put(0,0){\line(1,0){4}}
\end{picture}
\end{displaymath}
\bigskip
and $K$ be the \Fn\ diagram in (\ref{eq:DiagK}).
Then $$W_L(K)= -2W_L(P),$$

In other words
\be
W_L 
\DiagramK =
2W_L~\begin{picture}(4,1.5)
\qbezier(0.5,0)(2,2.5)(3.5,0)
\qbezier(1.5,0)(2,1)(2.5,0)
\thicklines
\put(0,0){\line(1,0){4}}
\end{picture}
-2W_L~
\begin{picture}(4,1.5)
\qbezier(0.5,0)(1.5,2)(2.5,0)
\qbezier(1.5,0)(2.5,2)(3.5,0)
\thicklines
\put(0,0){\line(1,0){4}}
\end{picture}%
~, 
\label{eq:lagrange}
\ee
i.e. 
\be 
\langle[a,b],[c,d]\rangle = 2\langle a,d \rangle \cdot \langle b,c
\rangle  -2 \langle a,c \rangle \langle b,d \rangle
\ \text{for} \ a,b,c,d \in sl_2 .
\label{eq:alglagr}
\ee

Under the isomorphism  $sl_2 \simeq so_3$, the inner product and the Lie
bracket on $sl_2$ become respectively the scalar and the vector product in
the three-dimensional space and~(\ref{eq:alglagr}) becomes the classical
Lagrange's identity%
\footnote{In fact, the identity~(\ref{eq:lagrange}) goes back to Euler's theory of
the motion of a rigid body. In the language of 20th century physics it looks like
{$\epsilon_{\a\b\g}\epsilon_{\rho\g\s} =\d_{\a\s}\d_{\b\rho} -
\d_{\a\rho}\d_{\b\s}, 
$}
where $\d$ is the Kronecker delta, and $\e$ is the standard completely antisymmetric tensor
in $\R^3$.}{}
$$\mathbf 
[a\times b]\cdot [c\times d] = (a\cdot c)(b\cdot d) - (a\cdot
d)(b\cdot c)
$$
which is equivalent to the better known {\em fundamental relation of vector calculus\/}
$$\mathbf
[a \times [b\times c]] = (a\cdot c) b - (a\cdot b) c~.
$$

Identity~(\ref{eq:lagrange}) allows us to compute the
universal $sl_2$ framed and unframed \ws s recursively.

\begin{thm} \label{thm:rec}
Let $D$ be a chord diagram, ``$a$'' a fixed chord in $D$, and
let $I_a = \{b_1,b_2,\ldots, b_p\}$ be 
the set of chords in $D$ intersecting $a$. Denote by
$D_a$ (resp.  $D_{a,i}$, $D_{a,ij}$) the diagram $D - a$ 
(resp. $D - a - b_i$, $D- a - b_i - b_j$). 
Then the framed $W$ and unframed
$\ov{W}$ universal $sl_2$ \ws s satisfy the following recursion relations

\begin{eqnarray}
(i) \quad  W(D) &=& (c-2p) \, W(D_a) + 2\sum_{i<j}
\left(W\left(D_{a,ij}^{||}\right) - 
W\left(D_{a,ij}^{\times}\right)\right)~,   \label{eq:recCV} \\
       & &\text{and}       \nonumber \\
\nonumber \\
(ii) \quad \ov{W}(D) &=&
-2\left| I_a\right| \ov{W}(D_a)
-2c\sum_{i\in I_a}\ov{W}(D_{a,i})
+2\sum_{i<j}
\left(
\ov{W}(D_{a,ij}^{||})-\ov{W}(D_{a,ij}^{\times})
\right)
\nonumber\\
 \label{eq:rec_defr} \\ 
& &
-2c \sum_{i<j}  \left(
\ov{W}(D_{a,ij}^{lr})+\ov{W}(D_{a,ij}^{rl})-\ov{W}(D_{a,ij}^{ll})-
\ov{W}(D_{a,ij}^{rr})
\right)~, 
\nonumber 
\end{eqnarray}
where $D_{a,ij}^{\times}$ (resp. $D_{a,ij}^{||}$)
 is the diagram obtained by adding to $D_{a,ij}$ two new
chords: the chord connecting the left end of $b_i$ with the right end of
$b_j$ and the chord connecting the left end of $b_j$ with the right end
of $b_i$ 
(resp. the cord connecting the left ends of $b_i$ and $b_j$ and the
chord connecting their right ends) 
 and  $D_{a,ij}^{lr}$ (resp. $D_{a,ij}^{rl}$, $D_{a,ij}^{ll}$, and
$D_{a,ij}^{rr}$) is the diagram obtained by adding to $D_{a,ij}$ a new
chord connecting the left end of $b_i$ and the right end of $b_j$
(resp. the right end of $b_i$ and the left end of $b_j$; the left ends
of $b_i$ and $b_j$; the right ends of $b_i$ and $b_j$) assuming that
the chord $a$ is drawn vertically. 
 (See Fig.~\ref{fig:diags} where 
$a$ is the vertical chord and $b_i,b_j$ are chords which
intersect $a$ so that $b_i$ is the upper chord and $b_j$ is the lower
chord.)
\begin{figure}[hbt]
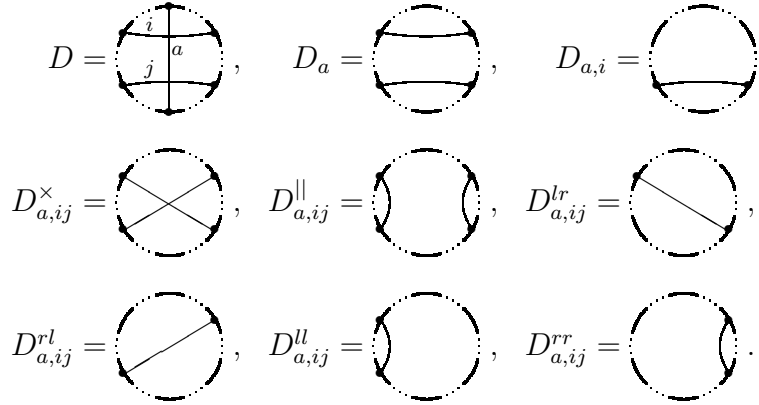


\begin{displaymath}
\begin{array}{rrr}
D=
\Picture{
\DottedCircle\FullChord[1,5]\FullChord[3,9]\FullChord[7,11]
\put(0.18,0.05){\makebox(0,0){$^a$}}
\put(-0.35,0.55){\makebox(0,0){$^i$}}
\put(-0.35,-0.25){\makebox(0,0){$^j$}}
}~,&
D_a= 
\Picture{
\DottedCircle\FullChord[1,5] \FullChord[7,11]\Arc[3]\Arc[9]
}~,&
D_{a,i} = \Picture{
\DottedCircle\FullChord[7,11]\Arc[3]\Arc[9]\Arc[1]\Arc[5]
} \\ [30pt]
D_{a,ij}^{\times} = \Picture{
\DottedCircle\FullChord[1,7]\FullChord[5,11]\Arc[3]\Arc[9]
}~,&
D_{a,ij}^{||} = \Picture{
\DottedCircle\FullChord[1,11]\FullChord[5,7]\Arc[3]\Arc[9]
}~,&
D_{a,ij}^{lr} = \Picture{
\DottedCircle\FullChord[5,11]\Arc[1] \Arc[3]\Arc[7]\Arc[9]
}~, \\ [30pt]
D_{a,ij}^{rl} = \Picture{
\DottedCircle\FullChord[1,7]\Arc[3]\Arc[5]\Arc[9]\Arc[11]
}~,&
D_{a,ij}^{ll} = \Picture{
\DottedCircle\FullChord[5,7]\Arc[1]\Arc[3]\Arc[9]\Arc[11]
}~,&
D_{a,ij}^{rr} = \Picture{
\DottedCircle\FullChord[1,11]\Arc[3] \Arc[5]\Arc[7]\Arc[9]
}~.
\end{array}
\end{displaymath}
\caption{Chord diagrams entering the recursion formulas}
\label{fig:diags}
\end{figure}
\end{thm}

Since all the diagrams in the right-hand side of~(\ref{eq:recCV})
and~(\ref{eq:rec_defr}) have 
fewer chords than $D$, this allows us to compute the value of the
\ws s $W$ and $\ov{W}$ on any chord diagram recursively. 

\begproof
Relation  (i)  for the framed invariant 
was proved by Chmutov and Varchenko in~\cite{CV}. 
It follows from  Lagrange's
identity~(\ref{eq:lagrange}) by induction on the number of chords
in $I_a$.

Let us prove relation (ii) for the unframed invariant.

By the deframing formula~(\ref{eq:deframed}) we have
\begin{eqnarray}
\ov{W}(D) & = & (U\cdot W)(D) = \sum_{E\subset D} U(D\setminus E)\cdot W(E) 
 =  \sum_{E\subset D} (-c)^{|D\setminus E|} W(E) \nonumber\\
  & = & \sum_{E \not \ni a}(-c)^{|D\setminus E|} W(E) 
      + \sum_{E \ni a}(-c)^{|D\setminus E|} W(E) \nonumber\\
  & = & \sum_{E\subset D_a}\left(
        (-c)^{|D\setminus E|} W(E) + (-c)^{|D\setminus E|-1} W(E+a) 
        \right)
\label{eq:ea} 
\end{eqnarray}

Applying  the recursion formula {\em (i)} to $W(E+a)$
in~(\ref{eq:ea}), we have 
\begin{eqnarray}
\ov{W}(D) & = &  \sum_{E\subset D_a}\biggl(
        (-c)^{|D\setminus E|} W(E) + (-c)^{|D\setminus E|-1} 
    \Bigl( (c-2|I_a\cap E|)W(E) \nonumber \\
&& \qquad  \qquad \qquad  \qquad \qquad  \qquad \qquad  \qquad 
+ \sum_{ {i<j,} \atop {b_i, b_j \in E} }
   2\bigl(W(E^{||}_{ij}) - W(E^{\times}_{ij})\bigr) \Bigr)
        \biggr)\nonumber\\
  & = &  \sum_{E\subset D_a}
-2(-c)^{|D\setminus E|-1}\biggl( |I_a\cap E|\cdot W(E) 
 - \sum_{ {i<j,} \atop {b_i,b_j\in E} }
    \left(W(E^{||}_{ij}) - W(E^{\times}_{ij})\right) \biggr)
\nonumber\\
  & = &  \sum_{E\subset D_a} \sum_{{i} \atop {b_i \in E}}
-2(-c)^{|D\setminus E|-1} \cdot W(E) 
 \label{eq:eb} \\
&& + \sum_{E\subset D_a}\sum_{ {i<j,} \atop {b_i,b_j\in E} }
2(-c)^{|D\setminus E|-1}
    \left(W(E^{||}_{ij}) - W(E^{\times}_{ij})\right)~. 
\label{eq:ec} 
\end{eqnarray}

Changing the order of summation in~(\ref{eq:eb}) and~(\ref{eq:ec}) we
get 
\begin{eqnarray}
\ov{W}(D) & = & \sum_{i} \sum_{{E\subset D_a} \atop {E\ni b_i}}
  -2 (-c)^{|D\setminus E|-1}W(E) \nonumber \\
&& + \sum_{i < j } \sum_{{E\subset D_a} \atop {E\ni b_i,b_j} }
2(-c)^{|D\setminus E|-1} 
\left( W(E^{||}_{ij}) - W(E^{\times}_{ij}) \right)\\
& = & \sum_{i} -2 S_i + \sum_{i<j}
2(S^{||}_{ij}-S^{\times}_{ij})~, \label{eq:sums1}
\end{eqnarray} 
where

\begin{eqnarray}
S_i &=& \sum_{{E\subset D_a} \atop {E\ni b_i} } 
   (-c)^{|D\setminus E|-1}W(E) \nonumber \\
&=&  \sum_{E\subset D_a} (-c)^{|D\setminus E|-1}W(E)  
- \sum_{E\subset D_{a,i} } (-c)^{|D\setminus E|-1}W(E)
\nonumber \\
&=&  \sum_{E\subset D_a} U(D_a\setminus E)W(E)  
    + \sum_{E\subset D_{a,i} } cU(D_{a,i}\setminus E)W(E),  
\nonumber 
\end{eqnarray}
and by the deframing formula (\ref{eq:deframed}) we obtain
\begin{eqnarray}
S_i &=& \ov{W}(D_a) +c\ov{W}(D_{a,i}).    \label{eq:term1}
\end{eqnarray}

For the term $S^{||}_{ij}$ in~(\ref{eq:sums1})
we have (by the inclusion-exclusion principle and 
the deframing formula)
\begin{eqnarray}
S^{||}_{ij} &=&   \sum_{{E\subset D_a} \atop {E\ni b_i,b_j} }
(-c)^{|D\setminus E|-1} W(E^{||}_{ij}) \nonumber \\
& = &  
\sum_{E\subset D^{||}_{a,ij}} (-c)^{|D\setminus E|-1} W(E) 
-   \sum_{E\subset D^{ll}_{a,ij}} (-c)^{|D\setminus E|-1} W(E)
                                                      \nonumber \\ 
& &  \qquad \qquad \qquad 
           -   \sum_{E\subset D^{rr}_{a,ij}}(-c)^{|D\setminus E|-1} W(E) 
                +   \sum_{E\subset D_{a,ij}} (-c)^{|D\setminus E|-1} W(E)
\nonumber \\
&=&
  \sum_{E\subset D^{||}_{a,ij}}  U(D^{||}_{a,ij}\setminus E)W(E)
+  \sum_{E\subset D^{ll}_{a,ij}} c U(D^{ll}_{a,ij}\setminus E)W(E)
         \nonumber \\
& & \qquad \qquad \qquad  
+  \sum_{E\subset D^{rr}_{a,ij}} c U(D^{rr}_{a,ij}\setminus E)W(E)
+  \sum_{E\subset D_{a,ij}} c^2 U(D_{a,ij}\setminus E)W(E)
\nonumber \\
&=&
\ov{W}( D^{||}_{a,ij}) +c\ov{W}( D^{ll}_{a,ij}) +c\ov{W}( D^{rr}_{a,ij})
+c^2 \ov{W}(D_{a,ij})~.                      \label{eq:term2}
\end{eqnarray}

Similarly, for $S^{\times}_{ij}$ in~(\ref{eq:sums1}) we get
\begin{eqnarray}
S^{\times}_{ij} &=& \sum_{{E\subset D_a} \atop {E\ni b_i,b_j} }
(-c)^{|D\setminus E|-1} W(E^{\times}_{ij}) \nonumber\\
&=&
  \ov{W}( D^{\times}_{a,ij}) - c\ov{W}(D^{ll}_{a,ij}) - c\ov{W}(D^{rr}_{a,ij})
- c^2 \ov{W}(D_{a,ij})    \label{eq:term3}
\end{eqnarray}

Substituting 
        (\ref{eq:term1}), (\ref{eq:term2}), and (\ref{eq:term3})
in (\ref{eq:sums1}),
we obtain the right hand side of~(ii).
\eproof

\begin{corol} \label{th:coef}
Let $\ov{w}_{j,n}$ be the order $n$ coefficients
of the universal \ws\ 
$$
   \ov{W}_{sl_2}=\sum\limits_{j,n} \ov{w}_{j,n}c^j
$$
 and
\be
\ov{W}_k = \sum_{n\ge 0} \ov{w}_{2n-k,n} \in \cW~.
\label{eq:lines}
\ee
Then 

(i) the \Vas\ series (\ref{eq:lines}) vanishes for $k < 0$, i.e.
\be
\ov{w}_{j,n}=0 \ \text{for} \  j>n/2~,
\label{eq:vanish}
\ee

(ii) the leading term series $\ov{W}_0=\sum\limits_{n\ge 0}
\ov{w}_{2n,n}$ 
satisfies the recursion relation
\be
\ov{W}_0(D)=-2\sum_i \ov{W}_0(D_{a,i})-2\sum_{i<j} 
\left( \ov{W}_0(D_{a,ij}^{lr}) + \ov{W}_0(D_{a,ij}^{rl}) -
\ov{W}_0(D_{a,ij}^{ll}) - \ov{W}_0(D_{a,ij}^{rr}) \right)~,
\label{eq:topline}
\ee

(iii) the  lower order series $\ov{W}_k$ satisfies the
recursion relation
\begin{eqnarray}
\ov{W}_k(D)&=&-2\sum_i \ov{W}_k(D_{a,i})-2\sum_{i<j} 
\left(  \ov{W}_k(D_{a,ij}^{lr}) + \ov{W}_k(D_{a,ij}^{rl}) - \ov{W}_k(D_{a,ij}^{ll} )-
\ov{W}_k(D_{a,ij}^{rr})\right)\nonumber\\
&& -2 |I_a|\ov{W}_{k-1}(D_a) +2\sum_{i<j} \left( \ov{W}_{k-1}(D_{a,ij}^{||})
- \ov{W}_{k-1}(D_{a,ij}^{\times}) \right)~. 
\end{eqnarray}
where the diagrams $D_a$, $D_{a,i}$,
 etc.\ are 
the same as in Theorem~\ref{thm:rec} (see Fig.~\ref{fig:diags}).

\end{corol}

\begproof
Relation (iii) follows immediately from~(\ref{eq:rec_defr})
and~(\ref{eq:lines}). Now (ii) and (i)  
follow from (iii) by induction on the number of chords in $D$, since
$\ov{W}_k(\Theta)=0$ for all $k$. 
\eproof

We will show that relations 
(\ref{eq:vanish}) and (\ref{eq:topline}) imply the \MM\ conjecture
(Theorem~\ref{thm:mm}), 
but first we have to recall some facts on the \AC\ \ws\ and its relation
with the universal \gl\ invariant.

\section{\AC\  and \gl\ universal \ws}

The universal \ws\ for the Lie superalgebra \gl\ was studied
in~\cite{FKV}. Here we recall the main results of~\cite{FKV}, 
and show that together with Corollary~\ref{th:coef} the
 recursion relation for $W_{\gl}$ implies the \MM\ conjecture.  

\subsection{The \gl\ \ws}

Let $V \cong \C^{1|1}$ be 
a $(1|1)$-dimensional  superspace.  
The Lie superalgebra of endomorphisms of $V$ is called
$\gl$. 

The bilinear form 
$$\langle x,y\rangle_{str} = str (xy)$$ 
on $\gl$ is invariant and
non-degenerate. 
  (The {\em supertrace\/} of a $(m|n)\times (m|n)$ matrix
$M = \pmatrix{A & B\cr C &  D\cr}$
is defined as   $ str M = tr A - tr B$).
Therefore, we can consider the universal \ws\
$W_{\gl}$ with values in $Z(\cU(\gl)) = \C[h,c]$,
where 
$h = \pmatrix{1 & 0\cr 0 & 1\cr} \in \gl$ 
and $c$ is the quadratic Casimir.

We have the following analogs of the fundamental
$sl_2$ relations~(\ref{eq:slkilling}) and~(\ref{eq:lagrange})
between invariant tensors on $L=\gl$ defined by the diagrams
$B$ and $K$ (see (\ref{eq:DiagK})). 

\be
W_L\Bubble = - 2 h^2.
\label{eq:gl1bb}
\ee
and
\be
W_L(K) = \half W_L(M),
\label{eq:gl11fund}
\ee
where
$$
M =
\begin{picture}(3,1)
\put(0.75,0){\oval(1,1)[tl]}
\put(1.25,0){\oval(1,1)[tr]}
\put(1,0.5){\circle{0.5}}
\put(2,0){\oval(1.5,1)[t]}
\thicklines
\put(0,0){\line(1,0){3}}
\end{picture} + 
\begin{picture}(3,1)
\put(1,0){\oval(1.5,1)[t]}
\put(1.75,0){\oval(1,1)[tl]}
\put(2.25,0){\oval(1,1)[tr]}
\put(2,0.5){\circle{0.5}}
\thicklines
\put(0,0){\line(1,0){3}}
\end{picture} - 
\begin{picture}(3,1.5)
\put(1.25,0){\oval(1,1)[tl]}
\put(1.75,0){\oval(1,1)[tr]}
\put(1.5,0.5){\circle{0.5}}
\put(1.5,0){\oval(2.5,2)[t]}
\thicklines
\put(0,0){\line(1,0){3}}
\end{picture} -
\begin{picture}(3,1.5)
\put(1.5,0){\oval(1.5,1)[t]}
\put(1.25,0){\oval(2,2)[tl]}
\put(1.75,0){\oval(2,2)[tr]}
\put(1.5,1){\circle{0.5}}
\thicklines
\put(0,0){\line(1,0){3}}
\end{picture}~.
$$
\medskip

Relation (\ref{eq:gl1bb}) shows in particular that the Killing form
for \gl\ is degenerate. The fundamental relation (\ref{eq:gl11fund})
was found in~\cite{FKV}, where we derived 
 the following recursive formulas for the universal $\gl$ \ws s.

\begin{thm}
Let $W$ be the universal \gl\ \ws.
In the notations of Theorem~\ref{thm:rec}
we have 

(i) $W(D)$ is a polynomial in $c$ and $h^2$ satisfying
\begin{eqnarray}
W(D) &=& c \, W(D_a) +h^2  \sum_i W(D_{a,i})    \label{eq:recgl} \\
&&{} -h^2 \sum_{i<j} \biggl(W\left(D_{a,ij}^{lr}\right) +
W\left(D_{a,ij}^{rl}\right)- W\left(D_{a,ij}^{ll}\right) -
W\left(D_{a,ij}^{rr}\right)\biggr)~, \nonumber
\end{eqnarray}
and

(ii)
the deframed \gl\ \ws\ $\ov{W}$ is given by setting $c=0$ in $W$
and satisfies
the following recursion relation.
\begin{equation}
\ov{W}(D)=h^2\sum_{i}\ov{W}(D_{a,i}) - h^2\sum_{i<j}\left( 
\ov{W}(D_{a,ij}^{lr})+\ov{W}(D_{a,ij}^{rl})-\ov{W}(D_{a,ij}^{ll})-%
\ov{W}(D_{a,ij}^{rr})\right). 
\label{eq:recac}
\end{equation}

\end{thm}
\medskip

The following result connects the \AC\ polynomial
 with the \gl\ universal invariant.

\begin{prop}[\cite{FKV}]
 The deframed \gl\ \ws\ coincides with the \Vas\ series
of the renormalized \AC\ invariant 
$\tilde{\D}$ (\ref{eq:alexnorm}).
\end{prop}

\subsection{\MM\ expansion}

Information on the universal $sl_2$ and \gl\ \ws s
allows us to study the canonical \VI\
appearing in the \MM\ expansion~(\ref{eq:expansion}).

\begin{prop}
The \Vas\ series $W_0=\sum\limits_{n\ge 0} w_{nn}$ of the diagonal
coefficients $v_{nn}$ of the \MM\ expansion~(\ref{eq:V0})
satisfies the recursion relation 
\be
W_0(D)=-\sum_i W_0(D_{a,i})-\sum_{i<j} 
\left( W_0(D_{a,ij}^{lr}) + W_0(D_{a,ij}^{rl}) -
W_0(D_{a,ij}^{ll}) - W_0(D_{a,ij}^{rr}) \right).
\label{eq:topjones}
\ee
\end{prop}
\begproof
The quadratic Casimir $c\in Z(\cU(sl_2))$ in the $d$-dimensional
irreducible representation acts as a multiplication by $(d^2-1)/2$.
Therefore, the \ws\ $w_{nn}$ can be expressed via the coefficients
$\ov{w}_{jn}$ of~(\ref{eq:lines}) as 
$$
w_{nn}=\cases{2^{-n}\ov{w}_{n,n/2}, &if $n$ is even,\cr
               0 & if $n$ is odd,}
$$
 and the
relation~(\ref{eq:topjones}) now follows  
from~(\ref{eq:topline}).
\eproof

Comparing the relations~(\ref{eq:topjones}) and~(\ref{eq:recac})
we come to a new proof of Proposition~\ref{prop:ws} 
and, therefore, of the  \MM\ conjecture.
\bigskip

\noindent
{\em Proof of Proposition~\ref{prop:ws}\/}.

Part (i) follows from Corollary~\ref{th:coef}.(i), 
since $\hW_{sl}(d) = \ov{W}_{sl_2}(\frac{d^2-1}{2})$
and the coefficients $w_{in}$ in (\ref{eq:coefmm})
  with $i > n$ are combinations of the coefficients
$\ov{w}_{jn}$ in~(\ref{eq:lines}) with $j > n/2$
and therefore vanish. 
\medskip

We will show that (\ref{eq:mm2}) follows from
recursion relations for $W_0$~(\ref{eq:topjones}) and
for $C$~(\ref{eq:recac}) by induction.

By the product formula~(\ref{eq:product}) and 
relations~(\ref{eq:topjones}) and~(\ref{eq:recac})
we have 
\begin{eqnarray*}
(W_0\cdot C)(D) &=& \sum_{E\subset D} W_0(D\setminus E)C(E)\\
&=& \sum_{E\subset D_a}\left( W_0(D \setminus E)C(E) + W_0(D_a
\setminus E) C(E+a) \right)\\
&=& \sum_{E\subset D_a} \biggl[  
\sum_{i \atop b_i\not \in  E} - W_0(D_{a,i}\setminus E) C(E) 
+ \sum_{i, \atop b_i \in E} W(D_a \setminus E)C(E_i)\\
&& -\sum_{i<j \atop b_i, b_j \not \in E}
\left(
  W_0(D_{a,ij}^{lr}\setminus E) + W_0(D_{a,ij}^{rl}\setminus E)   
 - W_0(D_{a,ij}^{ll}\setminus E)  - W_0(D_{a,ij}^{rr}\setminus E)  
\right) C(E)\\
&&-\sum_{i<j \atop b_i, b_j \in E} W_0(D_a \setminus E)
\left(  
C(E_{ij}^{lr}) + C(E_{ij}^{rl}) - C(E_{ij}^{ll}) - C(E_{ij}^{rr})
\right)
 \biggr]\\
&=& \sum_i \biggl(
\sum_{E\subset D_a \atop E \not \ni b_i}
-W_0(D_{a,i} \setminus E) C(E)
+ \sum_{E\subset D_a \atop E \ni  b_i}
W_0(D_a \setminus E) C(E_i)
\biggr)\\
&& -\sum_{i<j}
\biggl( \Bigl(
 \sum_{E \subset D_a \atop E \not \ni  b_i, b_j}
   W_0(D_{a,ij}^{lr}\setminus E)C(E)
+\sum_{E\subset D_a \atop E \ni b_i, b_j} 
       W_0(D_a \setminus E)C(E_{ij}^{lr}) \Bigr)\\
&&+ \Bigl(\sum_{E \subset D_a \atop E \not \ni  b_i, b_j}
   W_0(D_{a,ij}^{rl}\setminus E)C(E)
+\sum_{E\subset D_a \atop E \ni b_i, b_j} 
       W_0(D_a \setminus E)C(E_{ij}^{rl})\Bigr)\\
&&-\Bigl(\sum_{E \subset D_a \atop E \not \ni  b_i, b_j}
   W_0(D_{a,ij}^{ll}\setminus E)C(E)
+\sum_{E\subset D_a \atop E \ni b_i, b_j} 
       W_0(D_a \setminus E)C(E_{ij}^{ll})\Bigr)\\
&&-\Bigl(\sum_{E \subset D_a \atop E \not \ni  b_i, b_j}
   W_0(D_{a,ij}^{rr}\setminus E)C(E)
+\sum_{E\subset D_a \atop E \ni b_i, b_j} 
       W_0(D_a \setminus E)C(E_{ij}^{rr}) \Bigr) \biggr)\\
&=& 
   \sum_i\left(-(W_0\cdot C)(D_{a,i}) + (W_0\cdot C)(D_{a,i}) \right)\\
&& 
- \sum_{i<j}
 \left(
(W_0\cdot C)(D_{a,ij}^{lr}) + (W_0\cdot C)(D_{a,ij}^{rl})
- (W_0\cdot C)(D_{a,ij}^{ll}) - (W_0\cdot C)(D_{a,ij}^{rr})
  \right)\\
&=&
\sum_{i<j}
 \left(
(W_0\cdot C)(D_{a,ij}^{ll}) + (W_0\cdot C)(D_{a,ij}^{rr})
- (W_0\cdot C)(D_{a,ij}^{lr}) - (W_0\cdot C)(D_{a,ij}^{rl})
  \right)~.
\end{eqnarray*}

Now, since $(W_0\cdot C)(D)=1$ for $|D|=0$ and
$(W_0\cdot C)(D)=0$ for $|D|=1$
 we see by induction on $|D|$ that 
$(W_0\cdot C)(D)=0$ if $|D| \ge 1$.
\eproof

\end{document}